\documentclass[a4paper,12pt]{article}
\usepackage{amsfonts}
\usepackage[xdvi]{graphicx}
\usepackage[dvips]{color}
\usepackage{amsmath,amssymb,amsthm}

\newcommand{\R}{{\mathbb R}}
\newcommand{\C}{{\mathbb C}}

\makeatletter

\@addtoreset{equation}{section}

\newcounter{def}[section]
\renewcommand{\thedef}{\stepcounter{def}\thesection.\@arabic\c@def }

\begin{document}
\setlength{\baselineskip}{24pt}
\begin{center}
\textbf{\LARGE{A Hopf bifurcation in the Kuramoto-Daido model}}
\end{center}

\setlength{\baselineskip}{14pt}

\begin{center}
Institute of Mathematics for Industry, Kyushu University / JST PRESTO,\\
Fukuoka, 819-0395, Japan 
\end{center}
\begin{center}
Hayato CHIBA\footnote{chiba@imi.kyushu-u.ac.jp}
\end{center}
 \begin{center}
Oct 11, 2016
 \end{center}

\begin{center}
\textbf{Abstract}
\end{center}

A Hopf bifurcation in the Kuramoto-Daido model 
is investigated based on the generalized spectral theory and the center manifold reduction
for a certain class of frequency distributions.
The dynamical system of the order parameter on a four-dimensional center manifold is derived.
It is shown that the dynamical system undergoes a Hopf bifurcation as the coupling strength
increases, which proves the existence of a periodic two-cluster state of oscillators.


\section{Introduction}

\begin{figure}
\begin{center}
\includegraphics[]{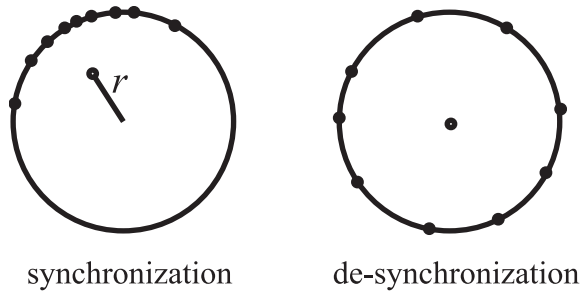}
\caption[]{Collective behavior of oscillators.}
\label{fig1}
\end{center}
\end{figure}

Collective synchronization phenomena are observed in a variety of areas such as chemical reactions,
engineering circuits and biological populations~\cite{Pik}.
In order to investigate such phenomena, a system of globally coupled phase oscillators called the
Kuramoto-Daido model \cite{Dai}
\begin{equation}
\frac{d\theta _i}{dt} 
= \omega _i + \frac{K}{N} \sum^N_{j=1} f (\theta _j - \theta _i),\,\, i= 1, \cdots  ,N,
\label{KM}
\end{equation}
is often used,
where $\theta _i = \theta _i(t) \in [ 0, 2\pi )$ is a dependent variable 
which denotes the phase of an $i$-th oscillator on a circle,
$\omega _i\in \R$ denotes its natural frequency drawn from some distribution function $g(\omega )$, $K>0$ is a coupling strength,
and where $f(\theta )$ is a $2\pi$-periodic function.
The complex order parameter defined by
\begin{eqnarray}
re^{i\psi} :=  \frac{1}{N}\sum^N_{j=1}e^{i\theta _j (t)}, \quad i = \sqrt{-1}
\label{1-2}
\end{eqnarray}
is used to measure the amount of collective behavior in the system;
if $r$ is nearly equal to zero, oscillators are uniformly distributed (called the incoherent state), while if $r>0$,
the synchronization occurs, see Fig. \ref{fig1}.

In this paper, the continuous limit (thermodynamics limit) of the following model
\begin{equation}
\frac{d\theta _i}{dt} 
= \omega _i + \frac{K}{N} \sum^N_{j=1}\bigl( \sin (\theta _j - \theta _i)
 + h\cdot \sin 2(\theta_j - \theta _i) \bigr) ,
\label{KMP}
\end{equation}
will be considered, where $h$ is a parameter which controls the strength of the second harmonic.
For the continuous limit of the system, a Hopf bifurcation from the incoherent state
to the two-cluster periodic state will be investigated based on the generalized spectral theory.

\begin{figure}
\begin{center}
\includegraphics[scale=1.3]{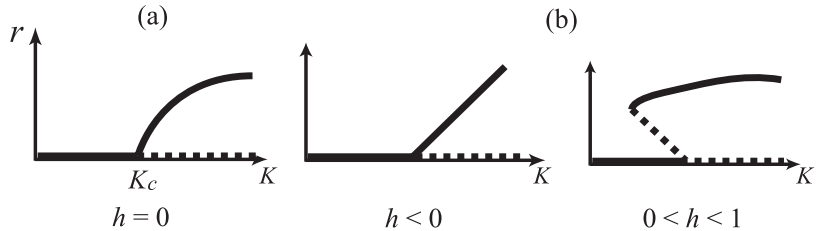}
\caption{Bifurcation diagrams of the order parameter for (a) $f(\theta ) = \sin \theta $
and (b) $f(\theta ) = \sin \theta  + h\sin 2\theta$.
The solid lines denote stable solutions, and the dotted lines denote unstable solutions. \label{fig2}}
\end{center}
\end{figure}

It is known that when the frequency distribution $g(\omega )$ is an even and unimodal function,
the transition from the incoherent state to the partially synchronized state 
occurs at the critical coupling strength $K=K_c = 2/(\pi g(0))$.
In Chiba \cite{Chi1, Chi2}, this result is proved based on the generalized spectral theory \cite{Chi3}
under the assumption that $g(\omega )$ has an analytic continuation near the real axis.
With the aid of the generalized spectral theory, it is proved that 
the order parameter is locally governed by the dynamical system on a center manifold as
\begin{eqnarray*}
\frac{dr}{dt} = \text{const.} \left(K-K_c + \frac{\pi g''(0) K_c^4}{16} r^2 \right) r+ O(r^4),
\end{eqnarray*}
for $h=0$, and
\begin{eqnarray*}
\frac{dr}{dt} = \text{const.} \left(K-K_c - \frac{K_c^2 Ch}{1-h}r \right) r + O(r^3),
\end{eqnarray*}
for $h\neq 0$, where $C$ is a certain negative constant.
As a result, a bifurcation diagram of $r$ is given as Fig.\ref{fig2}.
When $h=0$, the synchronous state emerges through a pitchfork bifurcation,
though when $h\neq 0$, it is a transcritical bifurcation.

The purpose in this paper is to investigate a Hopf bifurcation of the system (\ref{KMP})
under certain assumptions for the distribution function $g(\omega )$.
In particular, the dynamics of the order parameter on a center manifold will be derived.
For this purpose, we need five assumptions (A1) to (A5) given after Section 3.
Here, we give a rough explanation of these assumptions.
\\[0.2cm]
\textbf{(A1)} We assume that $h<1$ so that $\sin \theta $ is a dominant term in the coupling function.
\\
\textbf{(A2)} We assume that the distribution $g(\omega )$ of natural frequencies is an analytic function
near the real axis. This is the essential assumption to apply the generalized spectral theory.
\\
\textbf{(A3)} We will show that at a bifurcation value $K=K_c$, a pair of generalized eigenvalues 
of a certain linear operator obtained by the linearization of the system locates at 
the points $\pm iy_c$ on the imaginary axis.
We assume that such a pair is unique and they are simple eigenvalues.
\\
\textbf{(A4)} We assume that as $K$ increases, the pair of generalized eigenvalues 
transversally gets across the imaginary axis at the point $\pm iy_c$ from the left to the right.
\\
\textbf{(A5)} Assume that $g(\omega )$ is an even function.
\\

It seems that (A3) and (A4) are satisfied for a wide class of even and bimodal distributions
$g(\omega )$ as long as the distance of two peaks are sufficiently far apart, 
though we do not assume explicitly that $g(\omega )$ is bimodal.
The main results in the present paper are stated as follows;
\\[0.2cm]
\textbf{Theorem \thedef \,(Instability of the incoherent state).}
\\
Suppose (A1) and $g(\omega )$ is continuous.
There exists a number $\varepsilon >0$ such that
when $K_c<K<K_c+\varepsilon $, the incoherent state is linearly unstable,
where
\begin{eqnarray*}
K_c = \frac{2}{\pi g(y_c)}
\end{eqnarray*}
and $y_c$ is a certain real number, see (A3) above.
\\[0.2cm]
\textbf{Theorem \thedef\, (Local stability of the incoherent state).}
\\
Suppose (A1) and (A2).
When $0<K<K_c$, the incoherent state is linearly asymptotically stable in the weak sense
(see Section 4 for the weak stability).
\\[0.2cm]
\textbf{Theorem \thedef\, (Bifurcation).}
\\
Suppose (A1) to (A5).
There exists a positive constant $\varepsilon_0$ such that
if $K_{c}-\varepsilon _0 < K < K_{c} + \varepsilon_0 $ and if an initial condition is closed to the incoherent state, 
the dynamics of the order parameter is locally governed by a certain four dimensional dynamical system
on the center manifold given in Section 6.
At $K=K_c$ the system undergoes a Hopf bifurcation and when $K_{c}< K < K_{c} + \varepsilon_0,\, h\leq 0$
and $\mathrm{Re}(p_2) < 0$ (see below), the system has a family of asymptotically stable periodic orbits.

\textbf{(i)} Suppose $h=0$.
On the family of stable periodic orbits, the complex order parameter $\eta_1$ defined in Section 2 is given by
\begin{eqnarray*}
\eta_1(t) = 2 \sqrt{\frac{-\mathrm{Re}(p_1)}{\mathrm{Re}(p_2)}} \sqrt{K-K_c}\,e^{i\beta } \cos (y_c t + O(K-K_c)) + O(K-K_c),
\end{eqnarray*}
where $p_1$ and $p_2$ are certain complex constants explicitly given in Section 6,
and $\beta \in \R$ is an arbitrary constant specified by an initial condition.
The assumption (A4) implies $\mathrm{Re}(p_1) > 0$.

\textbf{(ii)} Suppose $h<0$.
On the family of stable periodic orbits, the complex order parameter $\eta_1$ is given by
\begin{eqnarray*}
\eta_1(t) = -2\frac{1-h}{hK_c}\mathrm{Re}(p_1)(K-K_c) \,e^{i\beta } \cos (y_c t + O(K-K_c)) + O((K-K_c)^2),
\end{eqnarray*}
where $p_1$ is the same constant as (i).
\\

The constants $p_1$ and $p_2$ are determined only by the frequency distribution $g(\omega )$.
The condition $\mathrm{Re}(p_2) < 0$ seems to be satisfied for most even and bimodal distributions.
If $\mathrm{Re}(p_2) > 0$, a family of unstable periodic orbits exists when
$K_{c}-\varepsilon _0 < K < K_{c}$; that is, a bifurcation occurs in the subcritical regime,
while the expression for $\eta_1$ is the same as above.
Similarly, if $0<h<1$, a bifurcation is subcritical and a family of unstable periodic orbits exists
for $K_{c}-\varepsilon _0 < K < K_{c}$.

In Martens et. al \cite{Mar}, the following bimodal frequency distribution defined as 
the sum of two Lorentzian distribution 
\begin{eqnarray}
g(\omega ) = \frac{1}{2\pi} \left( \frac{1}{(\omega -\omega _0)^2 +1} + \frac{1}{(\omega +\omega _0)^2 +1}\right),
\end{eqnarray}
is considered.
They revealed the dynamics of the order parameter in detail by using the Ott-Antonsen ansatz \cite{OA},
though it is applicable only when $h=0$.
For this bimodal distribution, we can verify that $K_c = 4,\, y_c = \sqrt{\omega _0^2 - 1}$,
$\mathrm{Re}(p_1) = 1/4$ and $\mathrm{Re}(p_2) = -4$.
Hence, there exists a family of stable periodic solutions for both of $h=0$ and $h<0$.
See also Example 3.5 and Example 5.3.


\section{The continuous model}

For the finite dimensional Kuramoto-Daido model (\ref{KM}), 
the $l$-th order parameter is defined by
\begin{eqnarray}
\hat{\eta}_l(t) := \frac{1}{N}\sum^N_{j=1} e^{i l \theta _j(t)}.
\label{2-1}
\end{eqnarray} 
By using it, Eq.(\ref{KM}) is rewritten as 
\begin{eqnarray*}
\frac{d\theta _j}{dt} = \omega _j + K \sum^\infty_{l=-\infty} f_l \hat{\eta}_l (t) e^{-i l \theta_j },
\quad f(\theta ) := \sum^\infty_{l=-\infty}f_l e^{il\theta }.
\end{eqnarray*}
The continuous limit of this model is an evolution equation of a density
\begin{equation}
\left\{ \begin{array}{ll}
\displaystyle \frac{\partial \rho_t}{\partial t} + \frac{\partial }{\partial \theta }(\rho_t v) = 0,
\quad \rho_t = \rho_t (\theta , \omega ), 
\label{conti}  \\
\displaystyle v := \omega  + K \sum^\infty_{l=-\infty} f_l \eta_l (t) e^{-i l \theta },  \\
\displaystyle \eta_l(t) := \int_{\R} \! \int^{2\pi}_{0} \!  e^{i l \theta } \rho_t (\theta , \omega ) g(\omega ) d\theta d\omega.
\end{array} \right.
\end{equation}
Here, $g(\omega )$ is a given probability density function for natural frequencies, 
and the unknown function $\rho_t = \rho_t (\theta , \omega )$ is a probability measure on $[0, 2\pi)$ parameterized by $t, \omega \in \R$.
$\eta_l(t)$ is a continuous analog of $\hat{\eta}_l(t)$ in (\ref{2-1}).
In particular, $\eta_1 (t)$ is a continuous version of Kuramoto's order parameter (\ref{1-2}).
The trivial solution $\rho_t = 1/(2\pi)$ of the system is a uniform distribution on the circle,
which is called the incoherent state (de-synchronous state).
Our purpose is to investigate the stability and bifurcation of the incoherent state
and the order parameter $\eta_1$.

Define the Fourier coefficients
\begin{eqnarray*}
Z_j(t, \omega ) := \int^{2\pi}_{0} \! e^{ij\theta }  \rho_t (\theta , \omega ) d\theta.
\end{eqnarray*}
Then, the continuous model is rewritten as a system of evolution equations of $Z_j$
\begin{eqnarray}
\frac{dZ_j}{dt} = ij\omega Z_j +i j K f_j \eta _j + i j K \sum_{l\neq j} f_l\eta_l Z_{j-l}.
\label{Z}
\end{eqnarray}
The trivial solution $Z_j \equiv 0 \, (j=\pm 1, \pm 2, \cdots)$ corresponds to the incoherent state
($Z_0 \equiv 1$ because of the normalization $\int^{2\pi}_{0} \! \rho_t (\theta , \omega ) d\theta =1$).
In what follows, we consider only the equations for $Z_1, Z_2,\cdots $ because $Z_{-j}$ is the complex
conjugate of $Z_j$.


\section{The transition point formula and linear instability}

To investigate the stability of the incoherent state, we consider the linearized system.
Let $L^2 (\R, g(\omega )d\omega )$ be the weighted Lebesgue space with the inner product
\begin{eqnarray*}
(\phi, \psi) = \int_{\R}\! \phi(\omega )\overline{\psi(\omega )} g(\omega )d\omega . 
\end{eqnarray*}
We define a one-dimensional integral operator $\mathcal{P}$ on $L^2 (\R, g(\omega )d\omega )$ by
\begin{equation}
(\mathcal{P}\phi ) (\omega ) = \int_{\R}\! \phi (\omega ) g(\omega )d\omega  = (\phi, P_0) \cdot P_0(\omega ),
\end{equation}
where $P_0 (\omega ) \equiv 1 \in L^2 (\R, g(\omega )d\omega )$ is a constant function.
Then, the order parameters are written by
\begin{equation}
\eta _j(t) 
= \int_{\R} \! Z_j(t, \omega )g(\omega ) d\omega = (Z_j, P_0)\cdot P_0(\omega ) = \mathcal{P}Z_j. 
\label{eta _j}
\end{equation}
Hence, Eq.(\ref{Z}) is expressed as
\begin{eqnarray}
\frac{dZ_j}{dt} = (ij\omega +i j K f_j \mathcal{P})Z_j + i j K \sum_{l\neq j} f_l (\mathcal{P}Z_l) Z_{j-l}.
\label{Z2}
\end{eqnarray}
The linearized system around the incoherent state is given by
\begin{equation}
\frac{dZ_j}{dt} = T_jZ_j := (ij\omega + ijKf_j \mathcal{P}) Z_j, \quad j=1,2,\cdots 
\label{linear}
\end{equation}
where $T_j = ij\omega + ijKf_j \mathcal{P}$ is a linear operator on $L^2(\R, g(\omega )d\omega )$.
Let us consider the spectra of $T_j$.
The multiplication operator $\phi (\omega ) \mapsto \omega \phi (\omega )$ on
$L^2(\R, g(\omega ) d\omega )$ is self-adjoint.
The spectrum of it consists only of the continuous spectrum given by $\sigma _c(\omega ) = \mathrm{supp} (g)$
(the support of $g$).
Therefore, the spectrum of the multiplication by $ij\omega $ lies on the imaginary axis;
$\sigma _c(ij\omega ) = ij\cdot \mathrm{supp} (g)$
(later we will suppose that $g$ is analytic, so that $\sigma _c(ij\omega )$ is the whole imaginary axis).
Since $\mathcal{P}$ is compact, it follows from the perturbation theory of linear operators \cite{Kato}
that the continuous spectrum of $T_j$ is given by $\sigma _c(T_j) = ij\cdot \mathrm{supp} (g)$,
and the residual spectrum of $T_j$ is empty.

When $f_j \neq 0$, eigenvalues $\lambda$ of $T_j$ are given as roots of the equation
\begin{equation}
\int_{\R} \! \frac{1}{\lambda - ij \omega }g(\omega )d\omega
= \frac{1}{ijKf_j}.
\label{eigen-eq}
\end{equation}
Indeed, the equation $(\lambda -T_j)v = 0$ provides
\begin{eqnarray*}
v+ijKf_j (v,P_0) (\lambda -ij\omega )^{-1} P_0 = 0.
\end{eqnarray*}
Taking the inner product with $P_0$, we obtain Eq.(\ref{eigen-eq}).
If $\lambda $ is an eigenvalue of $T_j$, the above equality shows that 
\begin{equation}
v_\lambda (\omega ) = \frac{1}{\lambda -ij\omega }
\label{ef}
\end{equation}
is the associated eigenfunction.
This is not in $L^2(\R, g(\omega )d\omega )$ when $\lambda $ is a purely imaginary number.
Thus, there are no eigenvalues on the imaginary axis.
Putting $\lambda =x+iy$ in Eq.(\ref{eigen-eq}) provides
\begin{eqnarray*}
\left\{ \begin{array}{ll}
\displaystyle \int_{\R}\! \frac{x}{x^2 + (y-j\omega )^2}g(\omega )d\omega
   = \frac{-\mathrm{Im}(f_j)}{jK|f_j|^2},  &  \\[0.4cm]
\displaystyle \int_{\R}\! \frac{y-j\omega }{x^2 + (y-j\omega )^2}g(\omega )d\omega 
  = \frac{\mathrm{Re}(f_j)}{jK|f_j|^2}, &  \\
\end{array} \right.
\end{eqnarray*}
which determines eigenvalues of $T_j$.
In what follows, we restrict our problem to the model (\ref{KMP}), for which
the coupling function is given by $f(\theta ) = \sin \theta + h\sin 2\theta $.
In this case, we have $f_1 = 1/(2i),\, f_2 = h/(2i)$ and $f_j = 0$ for $j\neq 1, 2$.
The spectrum of the operator $T_j$ for $j\neq 1, 2$ consists only of the continuous spectrum
on the imaginary axis.
$T_1$ and $T_2$ also have the continuous spectra on the imaginary axis.
Further, they have eigenvalues determined by the equations
\begin{equation}
\left\{ \begin{array}{ll}
\displaystyle \int_{\R}\! \frac{x}{x^2 + (y-\omega )^2}g(\omega )d\omega
   = \frac{2}{K},  & \\[0.4cm]
\displaystyle \int_{\R}\! \frac{y-\omega }{x^2 + (y-\omega )^2}g(\omega )d\omega = 0, &  \\
\end{array} \right.
\label{eigen1}
\end{equation}
and 
\begin{equation}
\left\{ \begin{array}{ll}
\displaystyle \int_{\R}\! \frac{x}{x^2 + (y-2\omega )^2}g(\omega )d\omega
   = \frac{h}{K},  & \\[0.4cm]
\displaystyle \int_{\R}\! \frac{y-2\omega }{x^2 + (y-2\omega )^2}g(\omega )d\omega = 0, &  \\
\end{array} \right.
\label{eigen2}
\end{equation}
respectively.
Eq.(\ref{eigen-eq}) for $j=1$ is given by
\begin{equation}
D(\lambda ) := \int_{\R}\! \frac{1}{\lambda -i\omega }g(\omega )d\omega = \frac{2}{K}.
\label{D} 
\end{equation}
The next lemma follows from formulae of the Poisson integral and the Hilbert transform.
\\[0.2cm]
\textbf{Lemma \thedef.} Suppose $g(\omega )$ is continuous. Then, the equality
\begin{eqnarray*}
\lim_{\lambda \to +0 + iy} D^{(n)}(\lambda ) 
&=& (-1)^n n! \cdot \lim_{\lambda \to +0+iy} \int_{\R}\! \frac{1}{(\lambda -i\omega )^{n+1}} g(\omega )d\omega \\
&=& \frac{1}{i^n} \cdot \lim_{\lambda \to +0+iy} \int_{\R}\! \frac{1}{\lambda -i\omega } g^{(n)}(\omega )d\omega \\
&=& \frac{1}{i^n}\left(  \pi g^{(n)} (y) - i\pi H[g^{(n)}](y)\right)
\end{eqnarray*}
holds for $n=0,1,2,\cdots $, where $\lambda \to +0+iy$ implies the limit to the point $iy \in i\R$
from the right half plane and $H[g]$ denotes the Hilbert transform defined by
\begin{eqnarray*}
H[g](y) &=& \frac{-1}{\pi} \mathrm{p.v.} \int_{\R} \frac{1}{\omega }g(\omega +y )d\omega \\
&=& \frac{-1}{\pi} \lim_{\varepsilon \to +0} \int^{\infty}_{\varepsilon }
\frac{1}{\omega }\left( g(y+\omega ) - g(y-\omega )\right) d\omega .
\end{eqnarray*}
\textbf{Lemma \thedef.} Suppose $K>0$. Then,
\\
(i) If an eigenvalue $\lambda $ of $T_1$ exists, it satisfies $\mathrm{Re} (\lambda ) > 0$.
\\
(ii) If $K>0$ is sufficiently large, there exists at least one eigenvalue $\lambda $ near infinity
on the right half plane.
\\
(iii) If $K>0$ is sufficiently small, there are no eigenvalues of $T_1$.
\\[0.2cm]
See \cite{Chi2, Chi4} for the proof.
\\

Eq.(\ref{eigen1}) combined with Lemma 3.1 yields
\begin{eqnarray*}
\left\{ \begin{array}{ll}
\displaystyle \lim_{x\to +0} \int_{\R}\! \frac{x}{x^2 + (y-\omega )^2}g(\omega )d\omega
    = \pi g(y)= \frac{2}{K},  & \\[0.4cm]
\displaystyle \lim_{x\to +0} \int_{\R}\! \frac{y-\omega }{x^2 + (y-\omega )^2}g(\omega )d\omega 
   = \pi H[g](y)= 0. &  \\
\end{array} \right.
\end{eqnarray*}
Let $y_1, y_2, \cdots $ be roots of the equation $H[g](y)= 0$, and put $K_j = 2/(\pi g(y_j))$.
The pair $(y_j, K_j)$ describes that some eigenvalue $\lambda =\lambda_j (K)$ of $T_1$
on the right half plane converges to the point $iy_j$ on the imaginary axis as $K \to K_j+0$.
Since $\mathrm{Re}(\lambda )>0$, the eigenvalue $\lambda_j (K)$ is absorbed into the 
continuous spectrum on the imaginary axis and disappears at $K=K_j$.
Suppose that $y_c$ satisfies $\sup_j \{ g(y_j)\} = g(y_c)$ and put
\begin{equation}
K_c = \inf_j \{ K_j\} = \frac{2}{\pi g(y_c)}.
\end{equation}
In what follows, $\lambda _c(K)$ denotes the eigenvalue of $T_1$ satisfying
$\lambda _c \to +0+iy_c$ as $K\to K_c+0$ ($y_c$ and $\lambda _c$ may not be unique).
The following formulae will be used later.
\\[0.2cm]
\textbf{Lemma \thedef.} The equalities
\begin{eqnarray*}
& & D(iy_c) := \lim_{\lambda \to +0+iy_c}D(\lambda ) = \frac{2}{K_c}, \\
& & \frac{d\lambda_c }{dK}\Bigl|_{K=K_c} = \frac{-2}{K_c^2 D'(iy_c)}
\end{eqnarray*}
hold.
\\[0.2cm]
\textit{Proof.} The first one follows from Eq.(\ref{D}) and the definition of $(y_c, K_c)$.
The derivative of Eq.(\ref{D}) as a function of $\lambda $ gives
\begin{eqnarray*}
D'(\lambda )= \frac{-2}{iK(\lambda )^2} \frac{dK}{d\lambda }.
\end{eqnarray*}
This proves the second one. $\Box$
\\

The eigenvalues of $T_2$ satisfy the same statement as Lemma 3.2.
The limit $x\to +0$ for Eq.(\ref{eigen2}) provides
\begin{eqnarray*}
\left\{ \begin{array}{ll}
\displaystyle \lim_{x\to +0}\int_{\R}\! \frac{x}{x^2 + (y-2\omega )^2}g(\omega )d\omega
   = \frac{1}{2}\pi g(y/2)= \frac{h}{K},  & \\[0.4cm]
\displaystyle \lim_{x\to +0}\int_{\R}\! \frac{y-2\omega }{x^2 + (y-2\omega )^2}g(\omega )d\omega 
   = \frac{1}{2}\pi H[g](y/2)= 0. &  \\
\end{array} \right.
\end{eqnarray*}
Let $y_1, y_2,\cdots $ be roots of the second equation, and
define $K_j^{(2)} = 2h /(\pi g(y_j/2))$ and $K_c^{(2)} = \inf_j \{ K_j^{(2)}\}$.
In what follows, we assume the following;
\\[0.2cm]
\textbf{(A1)} $h<1$.
\\

It is easy to verify that this condition is equivalent to $K_c < K_c^{(2)}$.
This implies that the eigenvalue $\lambda _c$ of $T_1$ still exists
on the right half plane after all eigenvalues of $T_2$ disappear as $K$ decreases.
In other words, as $K$ increases from zero, the eigenvalue $\lambda _c$ of $T_1$ first emerges
from the imaginary axis before some eigenvalue of $T_2$ emerges.
\\[0.2cm]
\textbf{Theorem \thedef\, (Instability of the incoherent state).}
\\
Suppose (A1) and $g(\omega )$ is continuous.
If $0<K<K_c$, the spectra of operators $T_1, T_2, \cdots $ consist only of the continuous spectra on the imaginary axis.
There exists a small number $\varepsilon >0$ such that
when $K_c<K<K_c+\varepsilon $, the eigenvalue $\lambda _c$ of $T_1$ exists on the right half plane.
Therefore, the incoherent state is linearly unstable.
\\

\begin{figure}
\begin{center}
\includegraphics[]{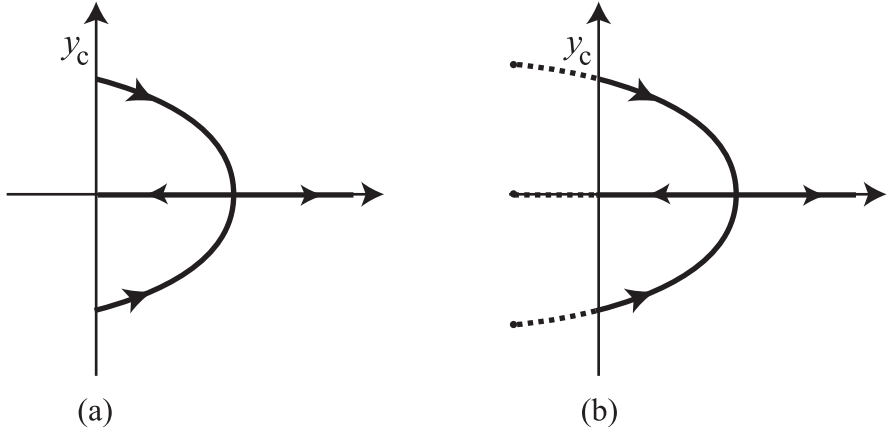}
\caption[]{ (a) The motion of the eigenvalues as $K$ increases for the distribution (\ref{L}).
The imaginary axis is the continuous spectrum. 
(b) The motion of the generalized eigenvalues as $K$ increases from zero for (\ref{L}).
The imaginary axis is a branch cut of the Riemann surface of the generalized resolvent.
The dotted curve denotes the path of the generalized eigenvalue on the second Riemann sheet.
See Section 5 for the detail.}
\label{fig3}
\end{center}
\end{figure}

This suggests that a first bifurcation occurs at $K=K_c$ and the eigenvalue $\lambda _c$ of $T_1$
plays an important role to the bifurcation.
\\[0.2cm]
\textbf{Example \thedef.} 
It is known that if $g(\omega )$ is an even and unimodal function, 
there exists a unique eigenvalue on the positive real axis for $K>K_c$.
Since we are interested in a Hopf bifurcation in this paper, 
let us consider the following bimodal frequency distribution defined as 
the sum of two Lorentzian distribution \cite{Mar}
\begin{eqnarray}
g(\omega ) = \frac{1}{2\pi} \left( \frac{1}{(\omega -\omega _0)^2 +1} + \frac{1}{(\omega +\omega _0)^2 +1}\right),
\label{L}
\end{eqnarray}
where $\omega _0>0$ is a parameter.
When $g''(0)>0 \Rightarrow \omega _0\geq 1/\sqrt{3}$, it is a bimodal function.
The equation $H[g](y) = 0$ has at most three roots given by
\begin{eqnarray*}
y_1 = 0, \quad y_2 = \sqrt{\omega _0^2 - 1}, \quad y_3 = -\sqrt{\omega _0^2 - 1}.
\end{eqnarray*}
Among them, $y_2$ and $y_3$ exist only when $\omega _0>1$.
Otherwise, the eigenvalue uniquely exists on the positive real axis as in the unimodal distribution case.
In what follows, we assume $\omega _0>1$.
Since $g(0)<g(y_2) = g(y_3) = 1/(2\pi)$, $y_c$ and $K_c$ are given by
\begin{eqnarray*}
y_c = \pm  \sqrt{\omega _0^2 - 1}, \quad K_c = \frac{2}{\pi g(y_c)} = 4.
\end{eqnarray*}
Eq.(\ref{D}) is calculated as
\begin{equation}
D(\lambda ) = \frac{\lambda +1}{(\lambda +1)^2 + \omega ^2_0} = \frac{2}{K}.
\label{ex}
\end{equation}
This shows that there are at most two eigenvalues on the right
half plane for each $K$.
The motion of the eigenvalues $\lambda =\lambda (K)$ as $K$ increases is represented in Fig.\ref{fig3} (a).
When $K<K_c = 4$, there are no eigenvalues.
At $K=K_c$, a pair of eigenvalues $\lambda _c(K_c) = \pm i \sqrt{\omega _0^2-1}$ pops up
from the continuous spectrum on the imaginary axis.
At $K=4\omega _0 > K_c$, two eigenvalues collide with one another on the real axis.
For $K>4\omega _0$, there are two eigenvalues on the positive real axis.
One of them goes to the left side as $K$ increases, and it is absorbed into the continuous 
spectrum and disappears at $K= 2/(\pi g(0))>4\omega _0$.
The other goes to infinity on the positive real axis as $K\to \infty$.
Later we will show that a Hopf bifurcation occurs at $K=K_c$.

\section{Linear stability}

When $0<K<K_c$, there are no spectra of operators $T_1, T_2, \cdots $ on the right half plane,
while the continuous spectra of them exist on the imaginary axis.
Hence, one may expect that the incoherent state is neutrally stable.
Nevertheless, we will show that the order parameter is asymptotically stable in a certain sense.
For this purpose, we need the following assumption.
Let $\delta $ be a positive number and define the stripe region on $\C$
\begin{eqnarray*}
S(\delta ) := \{ z\in \C \, | \, 0 \leq \mathrm{Im}(z) \leq \delta \}.
\end{eqnarray*}
We assume that
\\[0.2cm]
\textbf{(A2)} The distribution function $g(\omega )$ has an analytic continuation to the region $S(\delta )$.
On $S(\delta )$, there exists a constant $C>0$ such that the estimate
\begin{equation}
|g(z)| \leq \frac{C}{1+|z|^2}, \quad z\in S(\delta )
\end{equation}
holds.
\\

Let $H_+$ be the Hardy space on the upper half plane:
the set of bounded holomorphic functions on the real axis and the upper half plane.
It is a dense subspace of $L^2 (\R, g(\omega )d\omega )$.
For $\psi \in H_+$, set $\psi^*(z) := \overline{\psi (\overline{z})}$.

A function $f_t \in  L^2 (\R, g(\omega )d\omega )$ parameterized by $t$ is said to be convergent to zero
in the weak sense if the inner product $(f_t, \psi^*)$ decays to zero as $t\to \infty$ for any $\psi \in H_+$.
Note that $P_0\in H_+$ and the order parameter is written as $\eta_1(t) = (Z_1, P_0) = (Z_1, P_0^*)$.
This means that it is sufficient to consider the stability in the weak sense for the stability of the order parameter.
The next lemma plays an important role in the generalized spectral theory.
\\[0.2cm]
\textbf{Lemma \thedef.} Let $f(z)$ be a holomorphic function on the region $S(\delta )$.
Define a function $A[f](\lambda )$ of $\lambda $ to be 
\begin{eqnarray*}
A[f](\lambda ) = \int_{\R}\! \frac{1}{\lambda -i\omega }f(\omega )d\omega  
\end{eqnarray*} 
for $\mathrm{Re}(\lambda )>0$.
It has an analytic continuation $\hat{A}[f](\lambda )$ from the right half plane to the region 
$-\delta \leq \mathrm{Re}(\lambda ) \leq 0$ given by
\begin{equation}
\hat{A}[f](\lambda ) = \left\{ \begin{array}{ll}
A[f](\lambda ) & \mathrm{Re}(\lambda )>0 \\[0.2cm]
\displaystyle \lim_{\mathrm{Re}(\lambda ) \to +0}A[f](\lambda ) & \mathrm{Re}(\lambda ) = 0 \\[0.2cm]
A[f](\lambda ) + 2\pi f(-i\lambda ) & -\delta \leq \mathrm{Re}(\lambda ) < 0.
\end{array} \right.
\label{lemma4-1}
\end{equation}
See \cite{Chi2, Chi4} for the proof.
\\

\begin{figure}
\begin{center}
\includegraphics[scale=1.3]{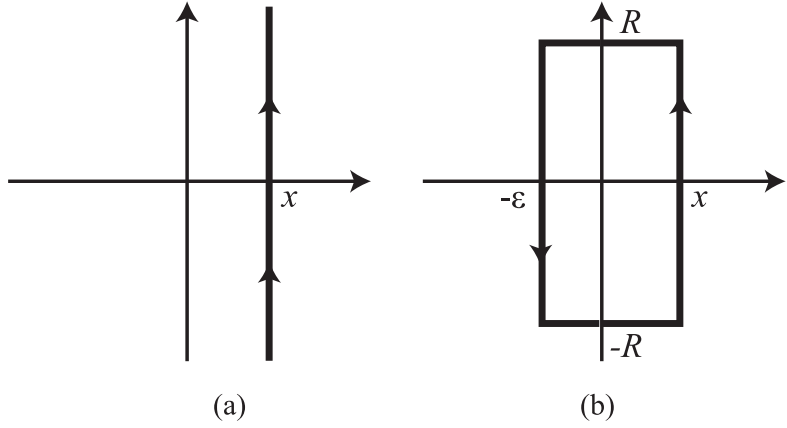}
\caption[]{Deformation of the integral path for the Laplace inversion formula.}
\label{fig4}
\end{center}
\end{figure}

It is known that the semigroup $e^{Tt}$ of an operator $T$ is expressed by the Laplace inversion formula
\begin{equation}
e^{Tt} = \lim_{y\to \infty} \frac{1}{2\pi i} \int^{x+iy}_{x-iy}\! e^{\lambda t}(\lambda -T)^{-1}d\lambda, 
\end{equation}
for $t>0$ (under a certain mild condition for $T$ \cite{Yos}).
Here, $x>0$ is chosen so that the integral path is to the right of the spectrum of $T$ (see Fig.\ref{fig4}(a)).
\\[0.2cm]
\textbf{Lemma \thedef.} The resolvent of $T_1= i\omega + iKf_1\mathcal{P}$ is given by
\begin{equation}
(\lambda -T_1)^{-1}\phi = (\lambda -i\omega )^{-1}\phi
 + \frac{iKf_1}{1-iKf_1D(\lambda )} ((\lambda -i\omega )^{-1}\phi, P_0) \frac{1}{\lambda -i\omega }.
\end{equation}
Let $\lambda _c$ be a simple eigenvalue of $T_1$.
The projection $\Pi_c$ to the eigenspace of $\lambda _c$ is given by
\begin{equation}
\Pi_c\phi = \frac{-1}{D'(\lambda _c)} ((\lambda _c - i\omega )^{-1} \phi, P_0) \frac{1}{\lambda _c - i\omega }.
\end{equation}
\\
See \cite{Chi4} for the proof.\\

Lemma 4.2 provides
\begin{eqnarray*}
& & ((\lambda -T_1)^{-1}\phi, \psi^*) \\
&=& ((\lambda -i\omega )^{-1}\phi, \psi^*)
 + \frac{iKf_1}{1-iKf_1D(\lambda )} ((\lambda -i\omega )^{-1}\phi, P_0)\cdot ((\lambda -i\omega )^{-1} \psi, P_0),
\end{eqnarray*}
which is meromorphic in $\lambda $ on the right half plane.
Suppose $\phi, \psi \in H_+$.
Due to Lemma 4.1, $((\lambda -T_1)^{-1}\phi, \psi^*)$
has an analytic continuation, possibly with new singularities, to the region 
$-\delta \leq \mathrm{Re}(\lambda ) \leq 0$
(Lemma 4.1 is applied to the factors $D(\lambda ), ((\lambda -i\omega )^{-1}\phi, \psi^*),
((\lambda -i\omega )^{-1}\phi, P_0)$ and $((\lambda -i\omega )^{-1}\psi, P_0)$).
A singularity on the left half plane is a root of the equation
\begin{equation}
1-iKf_1 (D(\lambda ) + 2\pi g(-i\lambda )) = 0.
\label{resonance}
\end{equation}
Such a singularity of the analytic continuation of the resolvent on the left half plane 
is called the generalized eigenvalue (see Sec.5 for the detail).

Now we can estimate the behavior of the semigroup by using the analytic continuation.
We have
\begin{eqnarray*}
(e^{T_1t}\phi, \psi^*) 
= \lim_{y\to \infty} \frac{1}{2\pi i} \int^{x+iy}_{x-iy}\! e^{\lambda t}((\lambda -T_1)^{-1}\phi, \psi^*) d\lambda, 
\end{eqnarray*}
where the integral path is given as in Fig.\ref{fig4} (a).
When $\phi, \psi \in H_+$, the integrand $((\lambda -T_1)^{-1}\phi, \psi^*)$ has an analytic continuation
to the region $-\delta \leq \mathrm{Re}(\lambda ) \leq 0$ which is denoted by $\mathcal{R}(\lambda )$.
\\[0.2cm]
\textbf{Lemma \thedef.}
Fix $K$ such that $0<K<K_c$.
Take positive numbers $\varepsilon, R$ and consider the rectangle shaped closed path $C$
represented in Fig.\ref{fig4} (b).
If $\varepsilon >0$ is sufficiently small, the analytic continuation of $((\lambda -T_1)^{-1}\phi, \psi^*)$
is holomorphic inside $C$ for any $R>0$.
\\
See \cite{Chi4} for the proof.\\

Because of this lemma, we have
\begin{eqnarray*}
0&=& \int^{x+iR}_{x-iR}\! e^{\lambda t} ((\lambda -T_1)^{-1}\phi, \psi^*) d\lambda
     + \int^{-\varepsilon -iR}_{-\varepsilon +iR} e^{\lambda t} \mathcal{R}(\lambda ) d\lambda \\
&+& \int^{iR}_{x+iR}\! e^{\lambda t} ((\lambda -T_1)^{-1}\phi, \psi^*) d\lambda
     + \int^{iR-\varepsilon }_{iR} e^{\lambda t} \mathcal{R}(\lambda ) d\lambda \\
&+& \int^{-iR+x}_{-iR}\! e^{\lambda t} ((\lambda -T_1)^{-1}\phi, \psi^*) d\lambda
     + \int^{-iR}_{-\varepsilon -iR} e^{\lambda t} \mathcal{R}(\lambda ) d\lambda .
\end{eqnarray*}
Due to the assumption (A2), we can verify that four integrals in the second and third lines above
become zero as $R\to \infty$.
Thus, we obtain
\begin{eqnarray*}
(e^{T_1t}\phi, \psi^*) 
= \lim_{R\to \infty} \frac{1}{2\pi i} \int^{-\varepsilon +iR}_{-\varepsilon -iR}\! 
   e^{\lambda t}\mathcal{R}(\lambda ) d\lambda.
\end{eqnarray*}
This proves $|(e^{T_1t}\phi, \psi^*) | \sim O(e^{-\varepsilon t})$ as $t\to \infty$.
We can show the same result for the operators $T_2, T_3,\cdots $.
\\[0.2cm]
\textbf{Theorem \thedef\, (Local stability of the incoherent state).}
\\
Suppose (A1) and (A2).
When $0<K<K_c$, $(e^{T_jt}\phi, \psi^*) $ decays to zero exponentially as $t\to \infty$
for any $j=1,2,\cdots $ and any $\phi, \psi\in H_+$.
Thus, the incoherent state is linearly asymptotically stable in the weak sense.

\section{The generalized spectral theory}

For the study of a bifurcation, we need generalized spectral theory developed in \cite{Chi3}
and applied to the Kuramoto model in \cite{Chi2} because the operator $T_1$
has the continuous spectrum on the imaginary axis (thus, the standard center manifold reduction is not applicable).
In this section, a simple review of the generalized spectral theory is given.
All proofs are included in \cite{Chi2, Chi3}.

Let $H_+$ be the Hardy space on the upper half plane with the norm
\begin{equation}
|| \phi ||_{H_+} = \sup_{\mathrm{Im}(z) > 0} |\phi (z)|.
\end{equation}
With this norm, $H_+$ is a Banach space.
Let $H_+'$ be the dual space of $H_+$; the set of continuous anti-linear functionals on $H_+$.
For $\mu \in H_+'$ and $\phi \in H_+$, $\mu (\phi )$ is denoted by $\langle \mu \,|\, \phi \rangle$.
For any $a,b \in \C,\, \phi, \psi \in H_+$ and $\mu, \xi \in H_+'$, the equalities
\begin{eqnarray*}
& & \langle \mu \,|\,  a \phi + b\psi\rangle 
   = \overline{a} \langle \mu \,|\,  \phi \rangle + \overline{b} \langle \mu \,|\, \psi \rangle, \\
& & \langle a\mu + b\xi \,|\, \phi \rangle
   = a \langle \mu \,|\, \phi \rangle + b \langle \xi \,|\, \phi \rangle,
\end{eqnarray*}
hold. An element of $H_+'$ is called a generalized function.
The space $H_+$ is a dense subspace of $L^2=L^2 (\R, g(\omega )d\omega )$
and the embedding $H_+ \hookrightarrow L^2$ is continuous.
Then, we can show that the dual $(L^2)'$ of $L^2$ is dense in $H_+'$ and it is continuously embedded in $H_+'$.
Since $L^2$ is a Hilbert space satisfying $ (L^2)' \simeq L^2$, we have
three topological vector spaces called a Gelfand triplet
\begin{eqnarray*}
H_+ \subset L^2 (\R, g(\omega )d\omega ) \subset H_+'.
\end{eqnarray*}
If an element $\phi \in H_+'$ is included in $L^2 (\R, g(\omega )d\omega )$, 
then $\langle \phi \,|\, \psi \rangle$ is given by 
\begin{eqnarray*}
\langle \phi \,|\, \psi \rangle := (\phi, \psi^*) = \int_{\R}\! \phi (\omega )\psi(\omega )g(\omega )d\omega . 
\end{eqnarray*}
(the conjugate $\psi^*$ is introduced to avoid the complex conjugate $\overline{\psi(\omega )}$
in the integrand).
Our operator $T_1$ and the above triplet satisfy all assumptions given in \cite{Chi3} to
develop a generalized spectral theory.
Now we give a brief review of the theory.
In what follows, we assume (A2).

The multiplication operator $\phi \mapsto i\omega \phi$ has the continuous spectrum 
on the imaginary axis; its resolvent is given by $(\lambda -i\omega )^{-1}$, and
it is not included in $L^2(\R, g(\omega )d\omega )$ 
when $\lambda $ is a purely imaginary number.
Nevertheless, we show that the resolvent has an analytic continuation from the right half plane to the 
left half plane in the generalized sense.
We define an operator $A(\lambda ) : H_+ \to H_+'$, parameterized by $\lambda \in \C$, to be
\begin{eqnarray*}
\langle A(\lambda )\phi \,|\, \psi \rangle
 = \left\{ \begin{array}{ll}
\displaystyle ((\lambda -i\omega )^{-1}\phi, \psi^*), & \mathrm{Re}(\lambda )>0, \\[0.4cm]
\displaystyle \displaystyle \lim_{\mathrm{Re}(\lambda ) \to +0}((\lambda -i\omega )^{-1}\phi, \psi^*)
   & \mathrm{Re}(\lambda )=0, \\[0.4cm]
\displaystyle ((\lambda -i\omega )^{-1}\phi, \psi^*) \\
\displaystyle  \quad +2\pi \phi (-i\lambda) \psi (-i\lambda)g(-i\lambda )& -\delta \leq \mathrm{Re}(\lambda )<0,
\end{array} \right.
\end{eqnarray*}
for $\phi, \psi \in H_+$.
Due to Lemma 4.1, $\langle A(\lambda )\phi \,|\, \psi \rangle$ is holomorphic.
That is, $A(\lambda )\phi$ is a $H_+'$-valued holomorphic function in $\lambda $.
In particular, $A(\lambda )$ coincides with $(\lambda -i\omega )^{-1}$ when $\mathrm{Re}(\lambda )>0$.
Since the continuous spectrum of the multiplication operator by $i\omega $ is the whole imaginary axis,
 $(\lambda -i\omega )^{-1}$ does not have an analytic continuation from the right half plane to the left half plane
as an operator on $L^2(\R, g(\omega )d\omega )$, however, 
it has a continuation $A(\lambda )$ if it is regarded as an operator from $H_+$ to $H_+'$.
$A(\lambda )$ is called the generalized resolvent of the multiplication operator by $i\omega $.

The next purpose is to define an analytic continuation of the resolvent of $T_1$ in the generalized sense.
Note that $(\lambda -T_1)^{-1}$ is rearranged as
\begin{eqnarray*}
(\lambda -i\omega -iKf_1\mathcal{P})^{-1}
= (\lambda -i\omega )^{-1}\circ (\mathrm{id} - iKf_1 \mathcal{P} (\lambda -i\omega )^{-1})^{-1}.
\end{eqnarray*}
Since the analytic continuation of $(\lambda -i\omega )^{-1}$ in the generalized sense is $A(\lambda )$,
we define the generalized resolvent $\mathcal{R}(\lambda ) : H_+ \to H_+'$ of $T_1$ by
\begin{eqnarray*}
\mathcal{R}(\lambda ):=A(\lambda ) \circ \left( \mathrm{id} - iKf_1 \mathcal{P}^\times A(\lambda ) \right)^{-1},
\end{eqnarray*}
where $\mathcal{P}^\times : H'_+ \to H'_+$ is the dual operator of $\mathcal{P}$.
For each $\phi \in H_+$, $\mathcal{R}(\lambda )\phi$ is a $H_+'$-valued meromorphic function.
It is easy to verify that when $\mathrm{Re}(\lambda )>0$, it is reduced to the usual resolvent $(\lambda -T_1)^{-1}$.
Thus, $\mathcal{R}(\lambda )$ gives a meromorphic continuation of $(\lambda -T_1)^{-1}$
from the right half plane to the left half plane as a $H_+'$-valued operator.
Again, note that $T_1$ has the continuous spectrum on the imaginary axis, so that it has no continuation
as an operator on $L^2(\R, g(\omega )d\omega )$.

A generalized eigenvalue is defined as a singularity of $\mathcal{R}(\lambda )$,
namely a singularity of $\left( \mathrm{id} - iKf_1 \mathcal{P}^\times A(\lambda ) \right)^{-1}$.
\\[0.2cm]
\textbf{Definition \thedef.}
If the equation
\begin{equation}
(\mathrm{id} - iKf_1 \mathcal{P}^\times A(\lambda ))\mu = 0
\end{equation}
has a nonzero solution $\mu$ in $H_+'$ for some $\lambda \in \C$, $\lambda $ is called a generalized eigenvalue
and $\mu$ is called a generalized eigenfunction.
\\

It is easy to verify that this equation is equivalent to
\begin{eqnarray}
\frac{2}{K} = \left\{ \begin{array}{ll}
D(\lambda ) & \mathrm{Re}(\lambda )>0,  \\
\displaystyle \lim_{\mathrm{Re}(\lambda ) \to +0}D(\lambda ) & \mathrm{Re}(\lambda )=0,  \\
D(\lambda )+2\pi g(-i\lambda ) & -\delta \leq \mathrm{Re}(\lambda )<0,
\end{array} \right.
\label{eigeneq3}
\end{eqnarray}
where we use $f_1 = 1/(2i)$.
When $\mathrm{Re}(\lambda ) > 0$, this is reduced to Eq.(\ref{D}).
In this case, $\mu$ is included in $L^2(\R, g(\omega )d\omega )$ and 
a generalized eigenvalue on the right half plane is an eigenvalue in the usual sense.
When $\mathrm{Re}(\lambda ) \leq 0$, this equation is equivalent to Eq.(\ref{resonance}).
The associated generalized eigenfunction is not included in $L^2 (\R, g(\omega )d\omega )$
but an element of the dual space $H'_+$.
Although a generalized eigenvalue is not a true eigenvalue of $T_1$, 
it is an eigenvalue of the dual operator:
\\[0.2cm]
\textbf{Theorem \thedef\, \cite{Chi2, Chi3}.}
Let $\lambda $ and $\mu$ be a generalized eigenvalue and the associated generalized eigenfunction.
The equality $T_1^\times \mu = \lambda \mu$ holds.
\\ 

Let $\lambda _0$ be a generalized eigenvalue of $T_1$ and 
$\gamma _0$ a small simple closed curve enclosing $\lambda _0$.
The generalized Riesz projection $\Pi_0 : H_+\to H_+'$ is defined by
\begin{eqnarray*}
\Pi_0 = \frac{1}{2\pi i}\int_{\gamma _0}\! \mathcal{R}(\lambda ) d\lambda .
\end{eqnarray*}
As in the usual spectral theory, the image of it gives the generalized eigenspace associated with $\lambda _0$.

\begin{figure}
\begin{center}
\includegraphics[scale=1.3]{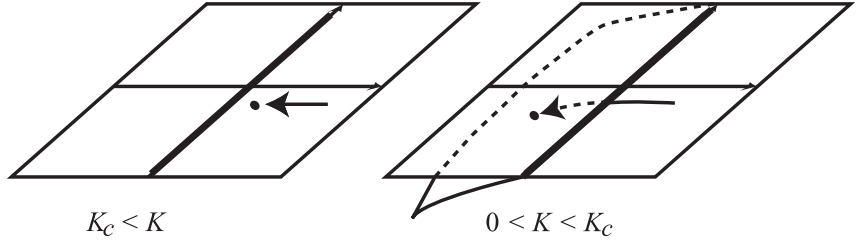}
\caption[]{The motion of the (generalized) eigenvalue as $K$ decreases.
When $0<K<K_c$, it lies on the second Riemann sheet of the resolvent and it is not a usual eigenvalue
but a generalized eigenvalue.}
\label{fig5}
\end{center}
\end{figure}

Let $\lambda = \lambda _c(K)$ be an eigenvalue of $T_1$ defined in Sec.3.
Recall that when $K_c<K$, $\lambda _c$ exists on the right half plane.
As $K$ decreases, $\lambda _c$ goes to the left side, and
at $K= K_c$, $\lambda _c$ is absorbed into the continuous spectrum on the imaginary axis and disappears.
However, we can show that even for $0<K<K_c$, 
$\lambda _c$ remains to exist as a root of Eq.(\ref{eigeneq3}) because the right hand side of Eq.(\ref{eigeneq3})
is holomorphic.
This means that although $\lambda _c$ disappears from the original complex plane at $K=K_c$,
it still exists for $0<K<K_c$ as a generalized eigenvalue on the Riemann surface of the generalized resolvent 
$\mathcal{R}(\lambda )$.
In the generalized spectral theory, the resolvent $(\lambda -T_1)^{-1}$ is regarded as an operator 
from $H_+$ to $H_+'$, not on $L^2(\R , g(\omega )d\omega )$.
Then, it has an analytic continuation from the right half plane to the left half plane
as $H_+'$-valued operator.
The continuous spectrum on the imaginary axis becomes a branch cut of the Riemann surface of the resolvent.
On the Riemann surface, the left half plane is two-sheeted (see Fig.\ref{fig5}).
We call a singularity of the generalized resolvent on the second Riemann sheet the generalized eigenvalue. 

On the dual space $H_+'$, the weak dual topology is equipped;
a sequence $\{ \mu_n \} \subset H_+'$ is said to be convergent to $\mu \in H_+'$
if $\langle \mu_n \,|\, \psi \rangle \in \C$ is convergent to $\langle \mu \,|\, \psi \rangle$ for each $\psi \in H_+$.
Recall that an eigenfunction of a usual eigenvalue $\lambda $ of $T_1$ is given by 
$v_\lambda (\omega ) = (\lambda -i\omega )^{-1}$ (Eq.(\ref{ef})).
A generalized eigenfunction $\mu_\lambda $ of a generalized eigenvalue $iy$ on the imaginary axis is given by
\begin{eqnarray*}
\mu_\lambda = \lim_{\lambda \to +0+iy} \frac{1}{\lambda -i\omega },
\end{eqnarray*}
where the limit is considered with respect to the weak dual topology.
This means that $\langle \mu_\lambda  \,|\, \psi \rangle$ is defined by
\begin{equation}
\langle \mu_\lambda  \,|\, \psi \rangle
 = \lim_{\lambda \to +0+iy} \langle \frac{1}{\lambda -i\omega } \,|\, \psi \rangle
 = \lim_{\lambda \to +0+iy} \int_{\R}\! \frac{1}{\lambda -i\omega }\psi (\omega )g(\omega )d\omega . 
\end{equation}
A generalized eigenfunction $\mu_\lambda $ associated with a generalized eigenvalue $\lambda $
on the left half plane is given by
\begin{equation}
\langle \mu_\lambda  \,|\, \psi \rangle
 =  \int_{\R}\! \frac{1}{\lambda -i\omega }\psi (\omega )g(\omega )d\omega + 2\pi \psi (-i\lambda )g(-i\lambda ). 
\end{equation}

To perform a center manifold reduction, we need the definition of a center subspace.
Usually, it is defined to be an eigenspace associated with eigenvalues on the imaginary axis.
For our case, the operators $T_1, T_2,\cdots $ have the continuous spectra on the imaginary axis.
Therefore, we define a generalized center subspace as a space spanned by generalized eigenfunctions
associated with generalized eigenvalues on the imaginary axis.
Note that this is a subspace of the dual $H_+'$, not of $L^2(\R, g(\omega )d\omega )$.
As $K$ increases from zero, some of the generalized eigenvalues of $T_1$ get across the imaginary axis
at $K=K_c$, and they become usual eigenvalues on the right half plane (see Fig.\ref{fig5}).
Hence, there is a nontrivial generalized center subspace at $K=K_c$ given by
\begin{eqnarray*}
\mathbf{E}^c := \mathrm{span}\{ \mu_\lambda \, | \, \lambda (K_c) \in i\R \},
\end{eqnarray*}
The next purpose is to perform a center manifold reduction.
\\[0.2cm]
\textbf{Example \thedef.} Let us consider the distribution (\ref{L}) given in Example 3.5.
The equation (\ref{eigeneq3}) for generalized eigenvalues is given by (\ref{ex});
the left hand side of it already gives an analytic continuation of $D(\lambda )$.
By solving it, it turns out that two generalized eigenvalues exist at $\lambda (0) = -1 \pm i\omega _0$
when $K=0$.
As $K$ increases, they go to the right side as is shown in Fig. \ref{fig3} (b).
They get across the imaginary axis when $K=K_c=4$, and become usual eigenvalues for $K>K_c$.
One of them again becomes a generalized eigenvalue at $K=2/(\pi g(0))$ by getting across
the imaginary axis from the right to the left.
The generalized center subspace for $K=K_c$ is a two-dimensional space.


\section{Center manifold reduction}

Recall that $y_c \in \R$ is defined as a number satisfying $\sup_j \{ g(y_j)\} = g(y_c)$,
where $y_1, y_2,\cdots $ are roots of the equation $H[g](y)=0$.
This gives a point $iy_c$ on the imaginary axis to which some eigenvalue of $T_1$ approaches as $K\to K_c+0$.
For a Hopf bifurcation, we assume the following:
\\[0.2cm]
\textbf{(A3)} There are exactly two nonzero values $y_c$ and $-y_c$ satisfying
$\sup_j \{ g(y_j)\} = g(\pm y_c)$.
Each of the corresponding eigenvalue of $T_1$ denoted by $\lambda _c^+(K)$ and $\lambda _c^-(K)$,
respectively, is simple near $K_c$ (i.e. the eigenspace is one dimensional).
\\[0.2cm]
\textbf{(A4)} The real part of $\displaystyle \frac{d\lambda_c^{\pm}}{dK}\Bigl|_{K=K_c}$ is positive.
\\[0.2cm]
\textbf{(A5)} $g(\omega )$ is an even function.
\\

The assumption (A3) implies that the generalized center subspace at $K=K_c$ is a two dimensional space given by 
\begin{eqnarray}
\mathbf{E}^c = \mathrm{span}\{ \mu_+, \mu_- \},
\quad \mu_{\pm} := \lim_{\lambda \to +0\pm iy_c}\frac{1}{\lambda -i\omega }.
\label{centersp}
\end{eqnarray}
The assumption (A4) means that the generalized eigenvalues $\lambda _c^{\pm}$ of $T_1$ transversely 
get across the imaginary axis from the left to the right.
Due to (A5), it is easy to verify that the following equalities hold:
\begin{equation}
D(iy_c) = D(-iy_c) = \frac{2}{K_c} ,\,\,
D'(iy_c) = \overline{D'(-iy_c)}, \,\, D''(iy_c) = \overline{D''(-iy_c)}.
\label{cc}
\end{equation}
It seems that (A3) and (A4) are satisfied for a wide class of even and bimodal distributions
$g(\omega )$ as long as the distance of two peaks are sufficiently far apart, see Example 3.5.
\\

In what follows, we assume (A1) to (A5).
We expect that a Hopf bifurcation occurs at $K=K_c$.
In Chiba \cite{Chi2}, the existence of the one dimensional center manifold in $H'_+$ is proved
for the Kuramoto model when $g(\omega )$ is even and unimodal.
In this paper, we formally perform the center manifold reduction without a proof of the 
existence of a center manifold.

We put $\varepsilon = K-K_c$, which plays a role of a bifurcation parameter.
Our ingredients are;
\\[0.2cm]
\textbf{Equations:} The equations (\ref{Z}) for $j=1,2$ with $f_1 = 1/(2i),\, f_2 = h/(2i)$ are given by
\begin{equation}
\left\{ \begin{array}{l}
\displaystyle \dot{Z}_1 = T_cZ_1 + \frac{\varepsilon }{2}\mathcal{P}Z_1 
 + \frac{K}{2} (h \cdot \eta_2 Z_{-1} -\overline{\eta_1} Z_2 -h\cdot \overline{\eta_2} Z_3), \\[0.2cm]
\displaystyle \dot{Z}_2 = T_2Z_2 + K (\eta_1Z_1 - \overline{\eta_1} Z_3 -h\cdot \overline{\eta_2} Z_4),
\end{array} \right.
\label{Y2}
\end{equation}
where $T_c$ is an operator $T_1$ estimated at $K = K_c$; that is,
$K$ in $T_1$ is denoted by $K = K_c + \varepsilon $ and accordingly $T_1 = T_c + \varepsilon \mathcal{P}/2$.
\\
\textbf{Center subspace:} As $K$ increases from zero, a pair of the generalized eigenvalues of 
$T_1$ denoted by $\lambda _c^{\pm} (K)$ gets across the imaginary axis at $\pm iy_c$ when $K=K_c$,
and they become usual eigenvalues on the right half plane when $K>K_c$.
The associated generalized eigenfunctions at $K=K_c$ and the generalized center subspace 
is given in (\ref{centersp}). 
\\
\textbf{Projection:} The projection to an eigenspace is given in Lemma 4.2.
The projection to the generalized center subspace spanned by $\mu_+$ and $\mu_-$ is
\begin{equation}
\Pi_c \phi = \frac{-1}{D'(iy_c)} \lim_{\lambda \to iy_c} ((\lambda  - i\omega )^{-1} \phi, P_0) \mu_+
+  \frac{-1}{D'(-iy_c)} \lim_{\lambda \to -iy_c} ((\lambda - i\omega )^{-1} \phi, P_0) \mu_-.
\label{6-4}
\end{equation}

We divide our result into two cases, $h=0$ and $h\neq 0$ because 
types of bifurcations of them are different.

\subsection{Center manifold reduction $(h=0)$}

Assume $h= 0$. Then, $T_2 = 2i\omega $.
Since $\Pi_cZ_1$ is a linear combination of $\mu_+$ and $\mu_-$,
we suppose $\Pi_cZ_1 = K_c/2 \cdot (\alpha_+ (t)\mu_+ + \alpha _-(t) \mu_- )$.
The scalar valued functions $\alpha_+ (t)$ and $\alpha _- (t)$ denote coordinates on the center subspace, 
and our purpose is to derive the dynamics of $\alpha_{\pm} $.
Since a solution decays to zero with an exponential rate for $(\mathrm{id}-\Pi_c)Z_1$ direction
and $Z_j\, (j=2,3,\cdots )$ directions, we assume that $(\mathrm{id}-\Pi_c)Z_1$ and $Z_j\, (j=2,3,\cdots )$
are of order $O( \alpha^2)$ which stand for $O(\alpha _+^2, \alpha _+\alpha _-, \alpha _-^2)$.
Thus, we write
\begin{equation}
Z_1 = \frac{K_c}{2} (\alpha_+ (t)\mu_+ + \alpha _-(t)\mu_- ) + O(\alpha ^2).
\label{6-5}
\end{equation}
Then, $\eta_1$ is given by
\begin{eqnarray}
\eta_1 &=& \int_{\R}\! Z_1 \cdot g(\omega )d\omega \nonumber \\
&=& \frac{K_c}{2}\alpha _+ \lim_{\lambda \to iy_c}\int_{\R}\! \frac{1}{\lambda -i\omega } g(\omega )d\omega 
    + \frac{K_c}{2}\alpha _- \lim_{\lambda \to -iy_c}\int_{\R}\! \frac{1}{\lambda -i\omega }g(\omega )d\omega 
+ O(\alpha ^2) \nonumber\\
 &=& \alpha _+ + \alpha _-  + O(\alpha ^2),
\label{6-6}
\end{eqnarray}
where we have used Eq.(\ref{cc}).
Further, we make the following ansatz
\begin{eqnarray}
\varepsilon \sim O(\alpha ^2) ,\quad \frac{d\alpha _{\pm}}{dt} = \pm iy_c \alpha _{\pm} + O(\alpha ^2),
\label{ansatz}
\end{eqnarray}
which will be verified if the dynamics on the center manifold is derived.

For $m,n = 0,1,2,\cdots $, we define functionals denoted by $\mu_+^m \cdot \mu_-^n \in H_+'$ by
\begin{equation}
\mu_+^m \cdot \mu_-^n := \lim_{\lambda _+ \to +iy_c} \lim_{\lambda _- \to -iy_c}
\frac{1}{(\lambda _+ - i\omega )^m} \frac{1}{(\lambda _- - i\omega )^n},
\end{equation}
where the limit is considered with respect to the weak dual topology.
\\[0.2cm]
\textbf{Lemma \thedef.} The following equalities hold.
\begin{eqnarray*}
& & \omega \mu_{\pm}^2 = i \mu_{\pm} \pm y_c \mu_{\pm}^2, \\
& & \omega \mu_{\pm} = i \pm y_c \mu_{\pm}, \\
& & \langle \mu_{\pm}P_0 \,|\, P_0 \rangle = D(\pm iy_c) = \frac{2}{K_c}, \\
& & \langle \mu_{\pm}^2 P_0 \,|\, P_0 \rangle = -D'(\pm iy_c), \\
& & \langle \mu_{\pm}^3 P_0 \,|\, P_0 \rangle = \frac{1}{2}D''(\pm iy_c), \\
& & \langle \mu_{+}\cdot \mu_- P_0 \,|\, P_0 \rangle = 0, \\
& & \langle \mu_{+}^2\cdot \mu_- P_0 \,|\, P_0 \rangle = \frac{1}{2iy_c} D'(iy_c), \\
& & \langle \mu_{+}\cdot \mu_-^2 P_0 \,|\, P_0 \rangle = \frac{-1}{2iy_c} D'(-iy_c),
\end{eqnarray*}
where $P_0(\omega ) = 1$ is a constant function.
\\[0.2cm]
\textit{Proof.} For the first equality, we have
\begin{eqnarray*}
\omega \mu_{\pm}^2 &=& \lim_{\lambda \to \pm iy_c}\frac{\omega }{(\lambda -i\omega )^2}
 = i \cdot \lim_{\lambda \to \pm iy_c}\frac{(\lambda -i\omega ) -\lambda }{(\lambda -i\omega )^2} \\
&=& i\cdot \lim_{\lambda \to \pm iy_c}\frac{1}{\lambda -i\omega } 
      \pm y_c \lim_{\lambda \to \pm iy_c}\frac{1}{(\lambda -i\omega )^2} \\
&=& i \mu_{\pm} \pm y_c \mu_{\pm}^2.
\end{eqnarray*}
The second one is proved in a similar manner.
The third one is given in Lemma 3.3.
The fourth and fifth equalities are easily shown by the integration by parts, see Lemma 3.1.
To prove the sixth equality, we use the partial fraction decomposition as
\begin{eqnarray*}
\langle \mu_{+}\cdot \mu_- P_0 \,|\, P_0 \rangle
&=& \lim_{\lambda _+ \to +iy_c} \lim_{\lambda _- \to -iy_c} \int_{\R}\! 
\frac{1}{(\lambda _+ - i\omega )(\lambda _--i\omega )}g(\omega )d\omega  \\
&=&  \lim_{\lambda _+ \to +iy_c} \lim_{\lambda _- \to -iy_c} \frac{1}{\lambda _+-\lambda _-}
\int_{\R}\! \left( \frac{-1}{\lambda _+-i\omega } + \frac{1}{\lambda _- - i\omega }\right) g(\omega )d\omega \\
&=& \frac{1}{2iy_c} \left( -D(iy_c) + D(-iy_c) \right) =0.
\end{eqnarray*}
The last two equalities are also verified by the partial fraction decomposition. $\Box$
\\[0.2cm]
\textbf{Lemma \thedef.}
Define 
\begin{eqnarray}
Z_2 = \frac{K_c^2}{4}\alpha _+^2 \mu_+^2 + \frac{K_c^2}{4}\alpha _-^2 \mu_-^2
      -\frac{K_c^2}{4iy_c}\alpha _+\alpha _- (\mu_+-\mu_-) + O(\alpha ^3).
\label{z2}
\end{eqnarray}
It satisfies the second differential equation of (\ref{Y2}) up to the order $O(\alpha ^3)$.
\\[0.2cm]
\textit{Proof.} By substituting Eqs.(\ref{6-5}), (\ref{6-6}) and (\ref{z2}) in the equation,
we can confirm with the aid of Lemma 6.1 and (\ref{ansatz}) that 
\begin{eqnarray*}
\dot{Z}_2 - (T_2Z_2 + K (\eta_1Z_1 - \overline{\eta_1} Z_3))
\end{eqnarray*}
is of order $O(\alpha ^3)$. $\Box$
\\

Let us apply the projection $\Pi_c$ to the both sides of the first equation of Eq.(\ref{Y2}) to get
\begin{eqnarray}
\frac{K_c}{2}(\dot{\alpha }_+ \mu_+ + \dot{\alpha }_- \mu_-)
= T_c^\times \Pi_c Z_1 + \frac{\varepsilon }{2} \eta_1 \Pi_c P_0 - \frac{K}{2}\overline{\eta}_1 \Pi_c Z_2. 
\label{6-9}
\end{eqnarray}
Theorem 5.2 gives
\begin{eqnarray*}
T_c^\times \Pi_c Z_1
= \frac{K_c}{2} T_c^\times (\alpha _+\mu_+ + \alpha _- \mu_-) 
= \frac{K_c}{2} \cdot iy_c \cdot (\alpha _+\mu_+ - \alpha _- \mu_-).
\end{eqnarray*}
The definition of $\Pi_c$ combined with Lemma 6.1 yields
\begin{eqnarray*}
\Pi_c P_0 &=& \frac{-2}{K_c D'(iy_c)} \mu_+ + \frac{-2}{K_c D'(-iy_c)} \mu_-, \\
\Pi_cZ_2 &=& \frac{-K_c^2}{4D'(iy_c)} \left( \frac{1}{2} \alpha _+^2 D''(iy_c)
 - \frac{1}{2iy_c}\alpha _-^2 D'(-iy_c) + \frac{1}{iy_c}\alpha _+\alpha _- D'(iy_c) \right) \mu_+ \\
& & + \frac{-K_c^2}{4D'(-iy_c)} \left( \frac{1}{2} \alpha _-^2 D''(-iy_c)
 + \frac{1}{2iy_c}\alpha _+^2 D'(iy_c) - \frac{1}{iy_c}\alpha _+\alpha _- D'(-iy_c) \right) \mu_- \\
& & + O(\alpha ^3). 
\end{eqnarray*}
Substituting these equalities into Eq.(\ref{6-9}) and comparing the coefficients of $\mu_+$ and $\mu_-$,
respectively, in the both sides of the equation, we obtain the dynamics on the center manifold 
\begin{equation}
\left\{ \begin{array}{ll}
\displaystyle \frac{d\alpha _+}{dt} = iy_c \alpha _+ + p_1 \varepsilon (\alpha _+ + \alpha _-) 
  + (\overline{\alpha }_+ + \overline{\alpha }_-) 
     \left( p_2\alpha _+^2 + p_3 \alpha _-^2 + p_4 \alpha _+\alpha _- \right) + O(\alpha ^4),  \\[0.4cm]
\displaystyle \frac{d\alpha _-}{dt} = -iy_c \alpha _- + \overline{p}_1 \varepsilon (\alpha _+ + \alpha _-) 
  + (\overline{\alpha }_+ + \overline{\alpha }_-) 
     \left( \overline{p}_2\alpha _-^2 + \overline{p}_3 \alpha _+^2 
      + \overline{p}_4 \alpha _+\alpha _- \right)+ O(\alpha ^4),  \\
\end{array} \right.
\label{centermfd}
\end{equation}
where $p_1$ to $p_4$ are complex numbers defined by
\begin{eqnarray*}
p_1 = \frac{-2}{K_c^2 D'(iy_c)}, \quad 
p_2 = \frac{K_c^2 D''(iy_c)}{8 D'(iy_c)}, \quad 
p_3 =-\frac{K_c^2 \overline{D'(iy_c)}}{8iy_cD'(iy_c)}, \quad 
p_4 = \frac{K_c^2}{4iy_c}.
\end{eqnarray*}
This is a (real) four dimensional dynamical system.
The next purpose is to reduce it.
Since the system is invariant under the action 
$(\alpha _+, \alpha _-) \mapsto (e^{i\beta }\alpha _+, e^{i\beta }\alpha _-)$ for $\beta \in \R$,
we can assume without loss of generality that $\mathrm{arg}(\alpha _+) + \mathrm{arg}(\alpha _-) = 0$.
Hence, we assume $\alpha _{\pm} = r_{\pm} e^{\pm i\psi}$ with $r_{\pm}, \psi \in \R$.
Substituting this into the system, we obtain the three dimensional system
\begin{equation}
\left\{ \begin{array}{l}
\dot{\psi} = y_c + O(r^2_{\pm}) = y_c + O(\varepsilon ),  \\[0.2cm]
\dot{r}_+ = \varepsilon \mathrm{Re}(p_1)r_+ + \varepsilon \mathrm{Re}(p_1 e^{-2i\psi}) r_-
 + \mathrm{Re}(p_2)r_+^3 + \mathrm{Re}(p_2e^{2i\psi})r_+^2r_- \\[0.2cm]
\quad +\mathrm{Re}(p_3e^{-4i\psi})r_+r_-^2 
 + \mathrm{Re}(p_3e^{-2i\psi})r_-^3 + \mathrm{Re}(p_4e^{-2i\psi})r_+^2r_- + O(r_{\pm}^4), \\[0.2cm]
\dot{r}_- = \varepsilon \mathrm{Re}(p_1)r_- + \varepsilon \mathrm{Re}(p_1 e^{-2i\psi}) r_+
 + \mathrm{Re}(p_2)r_-^3 + \mathrm{Re}(p_2e^{2i\psi})r_+r_-^2 \\[0.2cm]
\quad +\mathrm{Re}(p_3e^{-4i\psi})r_+^2r_- 
 + \mathrm{Re}(p_3e^{-2i\psi})r_+^3 + \mathrm{Re}(p_4e^{-2i\psi})r_+r_-^2 + O(r_{\pm}^4).
\end{array} \right.
\end{equation}
To derive this, note that $\mathrm{Re}(p_4) = 0$.
Now we apply the averaging method.
The right hand sides of the equations of $r_+$ and $r_-$ are averaged over $\psi$ to obtain
the averaging equation
\begin{equation}
\left\{ \begin{array}{l}
\dot{r}_+ = \varepsilon \mathrm{Re}(p_1)r_+ +  \mathrm{Re}(p_2)r_+^3 + O(r_{\pm}^4),  \\
\dot{r}_- = \varepsilon \mathrm{Re}(p_1)r_- +  \mathrm{Re}(p_2)r_-^3 + O(r_{\pm}^4). \\
\end{array} \right.
\label{averaging}
\end{equation}
It is known that the averaging equation provides an approximate solution within the error
of order $O(\varepsilon )$.
Further, if the averaging equation has a stable fixed point, 
then the original system has a stable periodic orbit \cite{Chi5}.
If $O(r^4_{\pm})$-terms are neglected, the averaging equation has at most four fixed points:
\begin{eqnarray*}
(r_+, r_-) = (0,0), \quad (r_*, 0), \quad (0,r_*), \quad (r_*, r_*),
\quad r_* := \sqrt{\frac{-\varepsilon \mathrm{Re}(p_1)}{\mathrm{Re}(p_2)}}.
\end{eqnarray*}
The last three fixed points exist as long as $-\varepsilon \mathrm{Re}(p_1)/\mathrm{Re}(p_2) > 0$.
The Jacobi matrices of the system at the fixed points are given by
\begin{eqnarray*}
\varepsilon \mathrm{Re}(p_1) \left(
\begin{array}{@{\,}cc@{\,}}
1 & 0 \\
0 & 1
\end{array}
\right),\,\, 
\varepsilon \mathrm{Re}(p_1) \left(
\begin{array}{@{\,}cc@{\,}}
-2 & 0 \\
0 & 1
\end{array}
\right),\,\, \varepsilon \mathrm{Re}(p_1) \left(
\begin{array}{@{\,}cc@{\,}}
1 & 0 \\
0 & -2
\end{array}
\right),\,\, \varepsilon \mathrm{Re}(p_1) \left(
\begin{array}{@{\,}cc@{\,}}
-2 & 0 \\
0 & -2
\end{array}
\right),
\end{eqnarray*}
respectively.
Because of the assumption (A4) and Lemma 3.3, we have $\mathrm{Re}(p_1) > 0$.
This shows that when $\varepsilon = K-K_c < 0$, the point $(0,0)$ is stable,
and when $\varepsilon =K-K_c>0$ and $\mathrm{Re}(p_2) < 0$, the fixed point $(r_*, r_*)$ exists and is stable.
This proves that when $\varepsilon>0$ and $\mathrm{Re}(p_2) < 0$, the averaging equation (\ref{averaging})
has a stable fixed point $(r_+, r_-) = (r_*, r_*) + O(\varepsilon )$,
and the system (\ref{centermfd}) has a family of stable periodic orbits
\begin{eqnarray*}
(\alpha _+, \alpha _-) = \left( r_* e^{i(y_ct + O(\varepsilon )) + i\beta }+ O(\varepsilon ),\,
r_* e^{-i(y_c + O(\varepsilon )) + i\beta } + O(\varepsilon ) \right),
\end{eqnarray*}
where $\beta\in \R$ is an arbitrary constant induced by the action 
$(\alpha _+, \alpha _-) \mapsto (e^{i\beta }\alpha _+, e^{i\beta }\alpha _-)$ and 
it is specified by an initial condition.
Since the order parameter is $\eta_1 = \alpha _+ + \alpha _- + O(\alpha ^2)$,
we obtain a family of stable solutions
\begin{equation}
\eta_1 = 2 r_* e^{i\beta } \cos (y_ct + O(\varepsilon )) + O(\varepsilon ).
\end{equation}
This completes the proof of Theorem 1.3 (i).

\subsection{Center manifold reduction $(h\neq 0)$}

Assume $h\neq 0$. Then, $T_2 = 2i\omega + hK\mathcal{P}$.
We again assume (\ref{6-5}), and make the following ansatz
\begin{eqnarray}
\varepsilon \sim O(\alpha ) ,\quad \frac{d\alpha _{\pm}}{dt} = \pm iy_c \alpha _{\pm} + O(\alpha ^2).
\label{ansatz2}
\end{eqnarray}
\\[0.2cm]
\textbf{Lemma \thedef.}
Define 
\begin{eqnarray}
Z_2 &=& \frac{K_c^2}{4}\alpha _+^2 \mu_+^2 + \frac{K_c^2}{4}\alpha _-^2 \mu_-^2
      -\frac{K_c^2}{4iy_c}\alpha _+\alpha _- (\mu_+-\mu_-) \nonumber \\
& & -\frac{hK_c^3}{8}\frac{D'(iy_c)}{1-h} \alpha _+^2 \mu_+ -\frac{hK_c^3}{8}\frac{D'(-iy_c)}{1-h} \alpha _-^2 \mu_-
+ O(\alpha ^3).
\label{z22}
\end{eqnarray}
It satisfies the second differential equation of (\ref{Y2}) up to the order $O(\alpha ^3)$.

This is proved in a similar manner to Lemma 6.2.
Let us apply the projection $\Pi_c$ to the both sides of the first equation of Eq.(\ref{Y2}).
\begin{eqnarray}
\frac{K_c}{2}(\dot{\alpha }_+ \mu_+ + \dot{\alpha }_- \mu_-)
= T_c^\times \Pi_c Z_1 + \frac{\varepsilon }{2} \eta_1 \Pi_c P_0 + \frac{K_ch}{2} \eta_2 \Pi_c Z_{-1} + O(\alpha ^3). 
\label{6-17}
\end{eqnarray}
Lemma 6.3 with Lemma 6.1 gives
\begin{eqnarray*}
\eta_2 = (Z_2, P_0) &=& -\frac{K_c^2}{4}\alpha _+^2 D'(iy_c) -\frac{K_c^2}{4}\alpha _-^2 D'(-iy_c)  \\
& & -\frac{hK_c^2}{4}\frac{D'(iy_c)}{1-h}\alpha _+^2 -\frac{hK_c^2}{4}\frac{D'(-iy_c)}{1-h}\alpha _-^2+O(\alpha ^3) \\
&=& -\frac{K_c^2}{4}\frac{D'(iy_c)}{1-h}\alpha _+^2 -\frac{K_c^2}{4}\frac{D'(-iy_c)}{1-h}\alpha _-^2+O(\alpha ^3). 
\end{eqnarray*}
\\[0.2cm]
\textbf{Lemma \thedef.} $\Pi_c Z_{-1}$ is given by
\begin{eqnarray*}
\Pi_cZ_{-1} = \frac{-4}{K_cD'(iy_c)} e^{-i \mathrm{arg}(\alpha _+)} \mu_+ 
  +  \frac{-4}{K_cD'(-iy_c)} e^{-i \mathrm{arg}(\alpha _-)} \mu_- + O(\alpha ). 
\end{eqnarray*}
See \cite{Chi1} for the proof.

Substituting these equalities into Eq.(\ref{6-17}) and comparing the coefficients of $\mu_+$ and $\mu_-$,
respectively, in the both sides of the equation, we obtain the dynamics on the center manifold 
\begin{equation}
\left\{ \begin{array}{ll}
\displaystyle \frac{d\alpha _+}{dt} = iy_c \alpha _+ + q_1 \varepsilon (\alpha _+ + \alpha _-) 
  + (q_2 \alpha _+^2 + q_3 \alpha _-^2) e^{-i \mathrm{arg}(\alpha _+)} + O(\alpha ^3),  \\[0.4cm]
\displaystyle \frac{d\alpha _-}{dt} = -iy_c \alpha _- + \overline{q}_1 \varepsilon (\alpha _+ + \alpha _-) 
  + (q_2 \alpha _-^2 + \overline{q}_3 \alpha _+^2) e^{-i \mathrm{arg}(\alpha _-)} + O(\alpha ^3) ,  \\
\end{array} \right.
\label{centermfd2}
\end{equation}
where $q_1, q_2$ and $q_3$ are complex numbers defined by
\begin{eqnarray*}
q_1 = \frac{-2}{K_c^2 D'(iy_c)}, \quad 
q_2 = \frac{hK_c}{1-h}, \quad 
q_3 =\frac{hK_c}{1-h}\frac{\overline{D'(iy_c)}}{D'(iy_c)}.
\end{eqnarray*}
($q_1$ is the same number as $p_1$). Note that $q_2$ is a real number.
The next purpose is to reduce this system by the same way as the last section.
Since the system is invariant under the action 
$(\alpha _+, \alpha _-) \mapsto (e^{i\beta }\alpha _+, e^{i\beta }\alpha _-)$ for $\beta \in \R$,
we can assume without loss of generality that $\mathrm{arg}(\alpha _+) + \mathrm{arg}(\alpha _-) = 0$.
Hence, we put $\alpha _{\pm} = r_{\pm} e^{\pm i\psi}$ with $r_{\pm}, \psi \in \R$.
Substituting this into the system, we obtain the three dimensional system
\begin{equation}
\left\{ \begin{array}{l}
\dot{\psi} = y_c + O(\varepsilon ),  \\[0.2cm]
\dot{r}_+ = \varepsilon \mathrm{Re}(q_1)r_+ + \varepsilon \mathrm{Re}(q_1 e^{-2i\psi}) r_-
 + q_2r_+^2 + \mathrm{Re}(q_3e^{-4i\psi})r_-^2 + O(r_{\pm}^3), \\[0.2cm]
\dot{r}_- = \varepsilon \mathrm{Re}(q_1)r_- + \varepsilon \mathrm{Re}(q_1 e^{-2i\psi}) r_+
 + q_2r_-^2 + \mathrm{Re}(q_3e^{-4i\psi})r_+^2 + O(r_{\pm}^3). \\[0.2cm]
\end{array} \right.
\end{equation}
Now we apply the averaging method.
The right hand sides of the equations of $r_+$ and $r_-$ are averaged over $\psi$ to obtain
the averaging equation
\begin{equation}
\left\{ \begin{array}{l}
\dot{r}_+ = \varepsilon \mathrm{Re}(q_1)r_+ + q_2r_+^2 + O(r_{\pm}^3),  \\
\dot{r}_- = \varepsilon \mathrm{Re}(q_1)r_- + q_2r_-^2 + O(r_{\pm}^3). \\
\end{array} \right.
\label{averaging2}
\end{equation}
If $O(r^3_{\pm})$-terms are neglected, the averaging equation has at most four fixed points:
\begin{eqnarray*}
(r_+, r_-) = (0,0), \quad (r_*, 0), \quad (0,r_*), \quad (r_*, r_*),
\quad r_* := \frac{-\varepsilon \mathrm{Re}(q_1)}{q_2}.
\end{eqnarray*}
The last three fixed points exist only when $-\varepsilon \mathrm{Re}(q_1)/q_2 > 0$.
The Jacobi matrices of the system at the fixed points are given by
\begin{eqnarray*}
\varepsilon \mathrm{Re}(q_1) \left(
\begin{array}{@{\,}cc@{\,}}
1 & 0 \\
0 & 1
\end{array}
\right),\,\, 
\varepsilon \mathrm{Re}(q_1) \left(
\begin{array}{@{\,}cc@{\,}}
-1 & 0 \\
0 & 1
\end{array}
\right),\,\, \varepsilon \mathrm{Re}(q_1) \left(
\begin{array}{@{\,}cc@{\,}}
1 & 0 \\
0 & -1
\end{array}
\right),\,\, \varepsilon \mathrm{Re}(q_1) \left(
\begin{array}{@{\,}cc@{\,}}
-1 & 0 \\
0 & -1
\end{array}
\right),
\end{eqnarray*}
respectively.
Because of the assumption (A4) and Lemma 3.3, we have $\mathrm{Re}(q_1) > 0$.
This shows that when $\varepsilon = K-K_c < 0$, the point $(0,0)$ is stable.
When $K-K_c<0$ and $q_2 > 0$, the fixed points $(r_*, 0),\, (0, r_*),\, (r_*, r_*)$ exist but they are unstable.
When $K-K_c>0$ and $q_2 < 0$, the fixed point $(r_*, r_*)$ is stable.
Since $h<1$ (the assumption (A1)), $q_2 < 0$ is equivalent to $h<0$.
This proves that when $K-K_c>0$ and $h < 0$, the averaging equation (\ref{averaging2})
has a stable fixed point $(r_+, r_-) = (r_*, r_*) + O(\varepsilon^2 )$,
and the system (\ref{centermfd2}) has a family of stable periodic orbits
\begin{eqnarray*}
(\alpha _+, \alpha _-) = \left( r_* e^{i(y_ct + O(\varepsilon )) + i\beta } +O(\varepsilon^2 ),\,
r_* e^{-i(y_c + O(\varepsilon )) + i\beta }+ O(\varepsilon^2 ) \right),
\end{eqnarray*}
where $\beta\in \R$ is an arbitrary constant induced by the action 
$(\alpha _+, \alpha _-) \mapsto (e^{i\beta }\alpha _+, e^{i\beta }\alpha _-)$ and 
it is specified by an initial condition.
Since the order parameter is $\eta_1 = \alpha _+ + \alpha _- + O(\alpha ^2)$,
we obtain a family of stable solutions
\begin{equation}
\eta_1 = 2 r_* e^{i\beta } \cos (y_ct + O(\varepsilon ))+ O(\varepsilon^2 ).
\end{equation}
This completes the proof of Theorem 1.3 (ii).



\begin{thebibliography}{99}
\setlength{\baselineskip}{0pt}

\bibitem{Chi1}
H.Chiba, I.Nishikawa,
Center manifold reduction for a large population of globally coupled phase oscillators, 
Chaos, 21, 043103 (2011).

\bibitem{Chi2}
H. Chiba,
A proof of the Kuramoto conjecture for a bifurcation structure of the infinite-dimensional Kuramoto model,
Ergodic Theory Dynam. Systems 35 (2015), no. 3, 762-834.

\bibitem{Chi3}
H. Chiba,
A spectral theory of linear operators on rigged Hilbert spaces under analyticity conditions, 
Adv. in Math. 273, 324-379, (2015).

\bibitem{Chi4}
H. Chiba,
A center manifold reduction of the Kuramoto-Daido model with a phase-lag,
(arXiv:1609.04126).

\bibitem{Chi5}
H. Chiba,
Extension and unification of singular perturbation methods for ODEs based on the renormalization group method, 
SIAM j. on Appl. Dyn.Syst., Vol.8, 1066-1115 (2009).

\bibitem{Dai}
H. Daido,
Onset of cooperative entrainment in limit-cycle oscillators with uniform all-to-all interactions:
bifurcation of the order function,
Phys. D 91, no. 1-2, 24-66, (1996).

\bibitem{Kato}
T. Kato,
Perturbation theory for linear operators,
Springer-Verlag, Berlin, 1995.

\bibitem{Mar}
E. A. Martens, E. Barreto, S. H. Strogatz, E. Ott, P. So, T. M. Antonsen,
Exact results for the Kuramoto model with a bimodal frequency distribution,
Phys. Rev. E (3) 79, 026204 (2009).

\bibitem{OA}
E. Ott, T. M. Antonsen,
Low dimensional behavior of large systems of globally coupled oscillators,
Chaos 18, 037113 (2008).

\bibitem{Pik}
A. Pikovsky, M. Rosenblum, J. Kurths, Synchronization: A Universal Concept
in Nonlinear Sciences, Cambridge University Press, Cambridge, 2001.

\bibitem{Yos}
K. Yosida,
Functional analysis,
Springer-Verlag, Berlin, 1995.


\end{thebibliography}
\end{document}